\documentstyle[12pt]{article}
\author{Youssef ALAOUI and My Abdelhakim EL IDRISSI SAAD}
\title{The Runge approximation theorem\\
for generalized polynomial hulls}
\date{}
\newcommand{\complexes}{\mbox{I}\!\!\!\mbox{C}}

\newtheorem{th}{theorem}
\newtheorem{lm}{lemma}
\begin{document}
\maketitle
\setcounter{page}{1}
{\large 1. Introduction}\\
\\
\hspace*{.1in}We consider for a compact set $K$ in $\complexes^{n}$
the polynomially convex hull
$$\hat{K}=\{z\in \complexes^{n}: |P(z)|\leq ||P||_{K} \ for \ all \ polynomials \ P\}$$
Where $||P||_{K}=sup\{|P(z)|:z\in K\}$.\\
\hspace*{.1in}It is known from the Runge approximation theorem that every
function which is holomorphic in a neighborhood of compact sets $K$
with $K=\hat{K}$ can be approximated uniformly on ${K}$ by analytic polynomials.\\
\hspace*{.1in}Our aim here is to prove the same result in the more general
situation when the r\^ole of $\hat{K}$ is played by the generalized polynomial
hull $h_{q}(K)$ introduced by Basener [1] and which can be defined, for each
integer $q\in \{0,1,...,n-1\}$, by
$$h_{q}(K)=\displaystyle\bigcap_{P\in \complexes[z_{1},\cdots ,z_{n}]}
\{z\in \complexes^{n}:|P(z)|\leq \delta_{K}(P,z)\}$$
Where $\delta_{K}(P,z)$ denotes the lowest value of $||P||_{K\cap f^{-1}(0)}$
when f ranges in the set of analytic polynomial maps $\complexes^{n}\rightarrow \complexes^{q}$
vanishing at z.\\
\hspace*{.1in}Notice, however, that this result does not give anything new when
$q=0$, since $h_{o}(K)=\hat{K}$.\newpage
\noindent{\large 2. Proof of the theorem}
\begin{lm}{- Let $K$ be a compact set in $\complexes^{n}$ with $K=h_{q}(K)$
and let $\omega$ be a neighborhood of $K$. Then one can find finitely
many holomorphic polynomials $P_{n+1},\cdots , P_{m}$ and a closed polydisc
$\Delta\subset \complexes^{n}$ such that
$$K\subset\displaystyle\bigcap_{j=n+1}^{m}\{z\in \Delta:|P_{j}(z)|\leq \delta_{K}(P_{j},z)\}=L\subset \omega$$}
\end{lm}
$Proof.$
\hspace*{.1in} We may, of course, assume that $\omega$ is bounded. Then there
exists a closed polydisc $\Delta(0,r)$ with center at $0$ and $r=(r_{1},
\cdots ,r_{n})$ such that $\omega\subset\Delta (0,r) $.
It follows that $$\Delta (0,r)\backslash\omega\subset\Delta (0,r)\backslash h_{q}(K)
\subset\displaystyle\bigcup_{P\in\complexes [z_{1},\cdots ,z_{n}]}
\{ z\in\complexes^{n}:P(z)>\delta_{K}(P,z)\}$$
We shall first prove that $\{ z\in \complexes^{n}:|P(z)|>\delta_{K}(P,z)\}$
is open. To see this, let $(z_{j})_{j\geq 1}$ be a sequence of points
in $K_{P}=\{z\in \complexes^{n}: |P(z)|\leq \delta_{K}(P, z)\}$ which converges
to a point $z\in \complexes^{n}$, and let $f: \complexes^{n}\longrightarrow \complexes^{q}$
be a polynomial map with $f(z)=0$. Define $f_{j}: \complexes^{n}\longrightarrow \complexes^{q}$
by $f_{j}=f-f(z_{j})$. Then\\ $|P(z_{j})|\leq ||P||_{K\cap f_{j}^{-1}(0)}$.
Let $\zeta_{j}$ be a point of $K\cap f_{j}^{-1}(0)$ such that\\
$|P(\zeta_{j})|=||P||_{K\cap f_{j}^{-1}(0)}$. Then there is subsequence
$\zeta_{j_{k}}$ converging towards a point $\zeta$. Since $|P(z_{j_{k}})|\leq |P(\zeta_{j_{k}})|$
and $f(\zeta_{j_{k}})\rightarrow 0$ when $k\rightarrow +\infty$,\\
then a passage to the limit shows that
$$|P(z)|\leq |P(\zeta)|\leq ||P||_{K\cap f^{-1}(0)}.$$
Hence $|P(z)|\leq \delta_{K}(P, z)$, and $z\in K_{P}$.\\
\hspace*{.1in}Because $\Delta (0,r)\backslash\omega$ is compact there are
$P_{n+1},\cdots ,P_{m}$ such that
$$\Delta\backslash\omega\subset \displaystyle\bigcup_{j=n+1}^{m}\{ z\in\complexes^{n}:
|P_{j}(z)|>\delta_{K}(P_{j},z)\}$$
Let $P_{j}(z)=\frac{z_{j}}{r_{j}}$ $j=1,\cdots ,n$ . Since $||P_{j}||_{K}\leq 1$ then
$$K\subset L=\displaystyle\bigcap_{j=n+1}^{m}\{ z\in\Delta:
|P_{j}(z)|\leq\delta_{K}(P_{j},z)\}\subset\omega$$\newpage
{\noindent \begin{lm}{- Let $L$ be a compact set in $\complexes^{n}$ and $f:\complexes^{n}
\longrightarrow\complexes^{q}$ a holomorphic polynomial map with $q$ in the
range $1\leq q\leq n-1$ such that $L\cap f^{-1}(0)\mp\emptyset$.\\ If
$\zeta\in \hat{L}\cap f^{-1}(0)$ is nonsingular point of $f^{-1}(0)$, then
$\zeta \in (L\cap f^{-1}(0))^{\hat{}}.$ \\Where $(L\cap f^{-1}(0))^{\hat{}}$
denotes the polynomially convex hull of the set $L\cap f^{-1}(0)$.}
\end{lm}}
{\noindent {\bf Proof.}\\
\\
\hspace*{.1in} Let $F:\complexes^{n}\times \complexes\longrightarrow\complexes^{q}$
be the map defined by $F(z,z_{n+1})=f(z)z_{n+1}$. By identifying $\complexes^{n}$
to $\complexes^{n}\times \{ 0\}$, $L$ can be considered as a compact subset of the analytic
variety $F^{-1}(0)=f^{-1}(0)\times \complexes \cup \complexes^{n}\times \{ 0\}$.
Since $\zeta$ lies in only one global branch $Z_{o}$ of $f^{-1}(0)$, then
$Z_{o}\times\complexes$ and $\complexes^{n}\times\{ 0\}$ are the unique
global branches of $F^{-1}(0)$ containing the point $(\zeta,0)$.
Set $Z=Z_{o}\times \complexes\cup \complexes^{n}\times \{0\}$.\\
If \ $(\zeta,0)\notin (L\cap Z)^{\hat{}}$ \ then there exists  \ a holomorphic
function $g$ on  $F^{-1}(0)$ such that $|g(\zeta ,0)|>1>||g||_{L\cap Z}$.
Let $h$ be a holomorphic function on $F^{-1}(0)$ such that $h(\zeta ,0)=1$ and $h=0 $ on $F^{-1}(0)\backslash Z$.
Then for sufficiently large positive integer $k$, the holomorphic function
$g_{k}=g^{k}h$ satisfies $|g_{k}(\zeta ,0)|>||g_{k}||_{L}$,\\ a contradiction.
We conclude that $(\zeta ,0)\in (L\cap Z)^{\hat{}}$.
\begin{lm}{- Let $K$ be a compact set in $\complexes^{n}$ with $K=h_{q}(K)$,
and $q$ in the range $1\leq q\leq n-1$. Let $\omega$ and $L$ be as in lemma 1.
Then $L$ is polynomially convex, ie. $L=\hat{L}$.}
\end{lm}
$Proof.$
\hspace*{.1in} The proof is by induction on $m-n$. Suppose at first that
$$L=\{ z\in \Delta:|P(z)|\leq\delta_{K}(P,z)\}$$
Given a point $z\in\hat{L}$ and an analytic polynomial map
$f:\complexes^{n}\longrightarrow \complexes^{q}$ with $f(z)=0$,
we shall prove that
$$|P(z)|\leq ||P||_{L\cap f^{-1}(0)}, \ if \ L\cap f^{-1}(0)\mp\emptyset$$
From lemma 2 we obtain the above inequality for any regular point\\
$z\in \hat{L}\cap f^{-1}(0)$ in $f^{-1}(0)$.\\
\hspace*{.1in}Now, let $z\in \hat{L}\cap Sing(f^{-1}(0))$, and assume that
$|P(z)|>||P||_{L\cap f^{-1}(0)}$.\\ Let $V$ be a small neighborhood of
$L\cap f^{-1}(0)$ such that $||P||_{V}<|P(z)|$.\\
By perturbing slightly the coefficients of the components of $f$, we can get
an analytic polynomial map $g=(g_{1},\cdots ,g_{q}):\complexes^{n}\longrightarrow\complexes^{q}$
with $g(z)=0$,\\ $z$ a regular point of $g^{-1}(0)$, and $g^{-1}(0)\cap L\subset V$.
(See [3]).\\ Which implies that $|P(z)|>||P||_{L\cap g^{-1}(0)}$, a contradiction.\\
Hence $|P(z)|\leq ||P||_{L\cap f^{-1}(0)}$ for all
$z\in \hat{L}\cap f^{-1}(0)$.\\
Let $\zeta $ be  a point of $L\cap f^{-1}(0)$ such that
$|P(\zeta )|=||P||_{L\cap f^{-1}(0)}$. Then, for all $z\in \hat{L}\cap f^{-1}(0),
|P(z)|\leq |P(\zeta)|\leq \delta_{K}(P,\zeta)\leq ||P||_{K\cap f^{-1}(0)}.$\\
Hence $|P(z)|\leq \delta_{K}(P,z)$, and $z\in L$.\\
\hspace*{.1in}We now assume that the statement of our lemma is already known for all positif integers
$\leq m-n-1$, and let $$L=\{ z\in L_{m-1}:|P_{m}(z)|\leq\delta_{K}(P_{m},z)\},$$
where
$L_{m-1}=\{z\in \Delta: |P_{j}(z)|\leq \delta_{K}(P_{j},z), j=n+1,\cdots ,m-1\}$.\\
Since, by the induction hypothesis, $L_{m-1}$ is polynomially convex, then
$$\hat{L}\subset \hat{L}_{m-1}\cap \{z\in \Delta: |P_{m}(z)|\leq \delta_{K}(P_{m},z)\}^{\hat{}}=L$$
\noindent\begin{th}{- Let $f$ be an analytic function in a neighborhood of a compact
set $K\subset\complexes^{n}$ with $K=h_{q}(K)$, and $q$ in the range $1\leq q\leq n-1$.
Then there is a sequence $f_{j}$ of analytic polynomials such that
$f_{j}\longrightarrow f$ uniformly on $K$}
\end{th}
$Proof.$
By lemma 3 we can find a compact polynomially convex set $L$\\ containing $K$
so that $f$ is analytic in a neighborhood of $L$. By the known version
of the Runge approximation theorem [2], $f$ can be approximated uniformly by
analytic polynomials on $L$. This completes the proof.\\
\\
\hspace*{.1in}{\small DEPARTEMENT DE MATHEMATIQUES, INSTITUT AGRONOMIQUE ET VETERINAIRE HASSAN II,
BP 6202 INSTITUTS (10101)\\ RABAT, MAROC}\\
\\
\hspace*{.1in}REFERENCES\\
\\
$[1]$ Basener, R.: Several dimensional properties of the spectrum of a\\ uniform
algebra. Pacific J. math. 74, 297-306 (1978).\\
$[2]$ Hormander, Lars.: An introduction to complex analysis in several \\variables.
Third, revised, edition, 1990.\\
$[3]$ Lupacciolu, G. and Stout.: Holomorphic hulls in termes of nonnegative
divisors and generalized polynomial hulls. Manuscripta Math. 98, 321-331 (1999)

\end{document}